\documentclass{interact}
\usepackage{epstopdf}
\usepackage{subfigure}
\usepackage{algorithm2e}%

\theoremstyle{plain}
\newtheorem{theorem}{Theorem}[section]
\newtheorem{corollary}[theorem]{Corollary}

\newtheorem{proposition}[theorem]{Proposition}
\newtheorem{alg}[theorem]{Algorithm}

\theoremstyle{definition}
\newtheorem{definition}{Definition}
\newtheorem{example}{Example}

\theoremstyle{remark}
\newtheorem{remark}{Remark}
\usepackage{enumitem}
\usepackage{color}
\begin{document}



\title{On the structure of join tensors with applications to tensor eigenvalue problems}

\author{
\name{Vesa Kaarnioja$^{\ast}$\thanks{$^\ast$Email: vesa.kaarnioja@aalto.fi}}
\affil{Aalto University School of Science, Department of Mathematics and Systems Analysis, FI-00076 Aalto, Finland}
}

\maketitle

\begin{abstract}
We investigate the structure of join tensors, which may be regarded as the multivariable extension of lattice-theoretic join matrices. Explicit formulae for a polyadic decomposition (i.e., a linear combination of rank-1 tensors) and a tensor-train decomposition of join tensors are derived on general join semilattices. We discuss conditions under which the obtained decompositions are optimal in rank, and examine numerically the storage complexity of the obtained decompositions for a class of LCM tensors as a special case of join tensors. In addition, we investigate numerically the sharpness of a theoretical upper bound on the tensor eigenvalues of LCM tensors. 
\end{abstract}


\begin{amscode}15A69; 06A12; 15B36 \end{amscode}

\section{Introduction}

Meet and join matrices are generalizations of classically studied GCD and LCM matrices, respectively, and they arise in lattice theory, where they are used to describe relations between distinct lattice elements. The study of GCD matrices began already in 1876, when Smith~\cite{smith} studied the determinant of the $n\times n$ matrix having the greatest common divisor of $i$ and $j$ as its~$ij$ element. The properties of meet and join matrices on semilattices have been studied in recent years by several authors: Ilmonen, Haukkanen, and Merikoski~\cite{i1} studied the eigenvalues of meet and join matrices associated with incidence functions. Mattila and Haukkanen~\cite{Mattila4} studied positive definiteness and eigenvalues of meet and join matrices, and Mattila, Haukkanen, and M\"{a}ntysalo~\cite{Mattila5} considered the singularity of LCM-type matrices.

Tensors (or multidimensional arrays) are a natural multivariable extension of matrices. Meet and join tensors can be used to describe relations between multiple lattice elements, and the attention of the lattice-theoretic community has been shifting toward tensors lately: Haukkanen~\cite{Ha3} considered the Cayley hyperdeterminants of GCD tensors and Luque~\cite{Luque} studied the Cayley hyperdeterminants of meet tensors.

Like meet and join matrices, the study of the properties of meet and join tensors relies heavily on explicit factorizations of these objects. For tensors, there exist many suitable candidates for representation in a factorized format. The canonical polyadic decomposition (also known as parallel factor analysis) is considered to be one of the fundamental tensor decompositions~\cite{KolBa09}. More recently, other formats have been proposed such as the tensor-train decomposition and the hierarchical Tucker decomposition (see, e.g.,~\cite{Osel11,Grasedyck10,Hackbusch14}). Tensor decompositions are useful for both theoretical and computational research work since they reveal the underlying structure of tensors arising from applications and, in consequence, can reduce the storage complexity of tensors drastically compared to the cardinality of the full array of tensor elements. Compressed representations of tensors are desirable from a numerical point of view since they remove the need to access all elements of the tensor array individually, thus enabling fast linear algebraic operations with tensors.

Recently, Ilmonen~\cite{ilhyp} discovered an explicit polyadic decomposition of meet tensors when the domain of definition is closed under the meet operation and generalized previously known bounds on the eigenvalues of meet matrices to formulate new bounds on the eigenvalues of meet tensors. The structure-theoretic result can be interpreted as a generalization of the $LDL^\textup{T}$ decompositions derived for meet and join tensors in the existing literature~\cite{Ov,Is2}. Subsequently in~\cite{hik16}, it was observed that factorizations of meet tensors on meet closed sets tend to be very sparse: for example, the storage of a class of $n$-dimensional order $d$ GCD tensors can be carried out in a low parametric format using only the order of $\mathcal{O}(n\ln n)$ elements, which is independent of the tensor order $d$ altogether. While the properties of meet tensors on meet closed sets have been demonstrated to be well suited for numerical computations involving high dimensionality and order, the case when the domain of definition is not closed has not been studied adequately. One example of a tensor that belongs to this category is the $n$-dimensional order $d$ LCM tensor
\[
A_{i_1,i_2,\ldots,i_d}={\rm lcm}(i_1,i_2,\ldots,i_d),\quad i_1,i_2,\ldots,i_d\in\{1,2,\ldots,n\},
\]
for which it was discovered in~\cite{hik16} that there exists (implicitly) a tensor-train decomposition with associated tensor-train ranks
\begin{align}
{\rm rank}_{\rm TT}A=\begin{cases}
n,&\text{if }d=3,\\
\#\{\alpha_1\cdots\alpha_{\lfloor d/2\rfloor}\mid 1\leq\alpha_i\leq n\text{ and }{\rm gcd}(\alpha_p,\alpha_q)=1,\ p\neq q\},&\text{if }d>3.
\end{cases}\label{lcmranktt}
\end{align}
The method of proof used to obtain this result is non-constructive and as such it does not provide an explicit tensor-train decomposition for LCM tensors.

The previous discussion lays the foundation for this paper, in which we derive new explicit representation formulae for join tensors in both polyadic and tensor-train formats without imposing the assumption that the domain of definition needs to be closed under the join operation. This approach covers, e.g., the aforementioned class of LCM tensors, which are a special case of join tensors defined on the divisor lattice. The approach taken in this paper generalizes the results of~\cite{ilhyp}, provides a constructive version of the result~\eqref{lcmranktt} for join tensors defined on general join semilattices,  and gives additional insight on the structure of join tensors. In this paper, the focus is kept on join tensors since it covers the important special case of LCM tensors in a natural way, but the results are straightforward to generalize to meet tensors as well.

In addition to deriving explicit decompositions of join tensors, we present a suffient condition which ensures that the obtained polyadic and tensor-train decompositions are rank decompositions, i.e., decompositions containing a minimal number of terms. We also conduct numerical experiments to examine the storage complexity of the obtained decompositions for a class of LCM tensors. 

Tensor eigenvalues have been an ongoing topic of research ever since their inception by Qi~\cite{Qi} and Lim~\cite{Lim}. The class of tensor eigenvalues studied by Qi and Lim inherit many properties from the theory of matrix eigenvalues. In particular, numerical schemes based on the power iteration~\cite{Ng09} and shifted power iteration~\cite{KoMa14} have been developed for the solution of extremal tensor eigenvalues. As an application of the decompositions developed in this paper, we consider tensor eigenvalues of join tensors and state a generalization of the bound discovered by Ilmonen~\cite{ilhyp} in the framework of general join tensors. We assess the sharpness of this upper bound for a class of LCM tensors in an ensemble of test cases utilizing the explicit tensor-train decomposition in the numerical solution of the dominant eigenvalues of LCM tensors.

This paper is organized as follows. In Section~\ref{prelims}, we overview the necessary definitions surrounding meet and join tensors and tensor decompositions. We move onto deriving explicit polyadic and tensor-train decompositions of join tensors in Section~\ref{joindecomp}. We investigate both the tensor-train rank and canonical rank of join tensors in Section~\ref{ttrank}. In Section~\ref{ttappl} numerical experiments are conducted on a class of LCM tensors to assess both their storage complexity and the sharpness of a theoretical upper bound. Finally, we end this paper on conclusions and discussion on possible future work in Section~\ref{concl}.

\section{Notations and preliminaries}\label{prelims}
\subsection{Semilattices}
Let $(P,\preceq)$ be a nonempty poset. The poset is \emph{locally finite} if the interval
\[
\{z\in P\mid x\preceq z\preceq y\}
\] 
is finite for all $x,y\in P$. If the greatest lower bound of $x,y\in P$ with respect to $\preceq$ exists, it is called the \emph{meet} of $x$ and $y$ and is denoted by $x\wedge y$. Conversely, if the least upper bound of $x,y\in P$ exists with respect to $\preceq$, it is called the \emph{join} of $x$ and $y$ and is denoted by $x\vee y$. If $x\wedge y\in P$ exists for all $x,y\in P$, then $(P,\preceq,\wedge)$ is called a \emph{meet semilattice}, and if $x\vee y\in P$ exists for all $x,y\in P$, then $(P,\preceq,\vee)$ is called a \emph{join semilattice}. If the poset $(P,\preceq,\wedge,\vee)$ is both a meet and a join semilattice, then it is called a \emph{lattice}.

Any function $f\!:P\times P\to\mathbb{C}$ such that $f(x,y)=0$, whenever $x\not\preceq y$, is called an \emph{incidence function} of $P$. If $f$ and $g$ are incidence functions of $P$, then their sum
\[
(f+g)(x,y)=f(x,y)+g(x,y),\quad x,y\in P,
\]
product
\[
(fg)(x,y)=f(x,y)g(x,y),\quad x,y\in P,
\]
and convolution
\[
(f*g)(x,y)=\sum_{x\preceq z\preceq y}f(x,z)g(z,y),\quad x,y\in P,
\]
are incidence functions of $P$ as well.

The incidence function $\delta$ defined by setting
\[
\delta(x,y)=\begin{cases}1,&\text{if }x=y,\\ 0&\text{otherwise},\end{cases}
\]
is unity under the convolution. The incidence function $\zeta$ of $P$ is defined by
\[
\zeta(x,y)=\begin{cases}1,&\text{if }x\preceq y,\\ 0&\text{otherwise.}\end{cases}
\]
The inverse of $\zeta$ under the convolution is called the M\"{o}bius function $\mu_P$ on $P$, and it can be computed inductively by
\begin{align*}
\mu_P(x,y)=\begin{cases}1,&\text{if }x=y,\\ -\sum_{x\preceq z\prec y}\mu_P(x,z),&\text{if }x\prec y,\\
0&\text{otherwise}\end{cases}
\end{align*}
for all $x,y\in P$.

For further material on arithmetic functions and incidence algebras, the reader may be interested to see~\cite{Ai,Mc,St}.

\subsection{Tensor decompositions}
An $n$-dimensional order $d$ tensor $A$ is defined as an array of $n^d$ elements
\[
A_{i_1,i_2,\ldots,i_d},
\]
where $i_1,i_2,\ldots,i_d\in\{1,2,\ldots,n\}$. The tensor is said to be symmetric if its entries are invariant under any permutation of its indices.

Let $x=(x_1,x_2,\ldots,x_n)^\textup{T}$. The mode-$k$ contraction of the tensor $A$ against the vector $x$ is defined as the tensor $B=A\times_k x$ defined by setting
\[
B_{i_1,\ldots,i_{k-1},1,i_{k+1},\ldots,i_d}=\sum_{i_k=1}^nA_{i_1,\ldots,i_k,\ldots,i_d}x_{i_k},\quad i_1,\ldots,i_{k-1},i_{k+1},\ldots,i_d\in\{1,\ldots,n\}.
\]
In addition, we define the scalar $Ax^d$ and $n$-vector $Ax^{d-1}$, respectively, by setting
\[
Ax^d=A\times_1x\times_2\cdots\times_dx\quad\text{and}\quad Ax^{d-1}=A\times_2x\times_3\cdots\times_dx.
\]
For $\alpha>0$, we denote $x^{[\alpha]}=(x_1^\alpha,x_2^\alpha,\ldots,x_n^\alpha)^\textup{T}$.

Unlike the case of vectors and matrices, which are simply first and second order tensors, respectively, it is inefficient to consider higher order tensors as simple arrays of values since this would require storing a total of $n^d$ elements. Because of the exponential dependence on $d$, the number of elements that need to be stored is prohibitively high for even moderate order tensors.

In the case of symmetric tensors, it is sufficient to store only the so-called symmetric part. The matrix analogue of the symmetric part is equivalent to storing only the lower triangular part of a symmetric matrix.
\begin{definition}[Symmetric part]
Let $A$ be a symmetric $n$-dimensional order $d$ tensor. Then its symmetric part is defined by setting
\[
{\rm Sym}\, A=\sum_{i_1=1}^n\sum_{i_2=1}^{i_1}\cdots\sum_{i_d=1}^{i_{d-1}}A_{i_1,\ldots,i_d}{\rm e}^{i_1}\otimes\cdots\otimes {\rm e}^{i_d},
\]
where ${\rm e}^{i}=[\delta(1,i),\ldots,\delta(n,i)]^\textup{T}$ denotes the $i^{\rm th}$ Euclidean unit vector and the Segre outer product of vectors $v^1,\ldots,v^d\in\mathbb{C}^n$ is defined as the $n$-dimensional order $d$ tensor defined elementwise by setting $(v^1\otimes\cdots\otimes v^d)_{i_1,\ldots,i_d}=v_{i_1}^1\cdots v_{i_d}^d$, $i_1,\ldots,i_d\in\{1,\ldots,n\}$.
\end{definition}
\begin{remark}\label{symremark}
Note that ${\rm Sym}\, A\neq A$. The symmetric part can be used to recover the full tensor via
\[
A_{i_{\sigma(1)},\ldots,i_{\sigma(d)}}=({\rm Sym}\,A)_{i_1,\ldots,i_d}\quad\text{for }1\leq i_1\leq\cdots\leq i_d\leq n,
\]
where $\sigma$ is any permutation of $\{1,2,\ldots,d\}$.
\end{remark}

Even accounting for symmetry, the storage of
\[
\#\{(i_1,\ldots,i_d)\mid 1\leq i_1\leq\cdots\leq i_d\leq n\}=\binom{d+n-1}{d}=\frac{1}{d!}\prod_{k=0}^{d-1}(n+k)=\mathcal{O}(n^d/d!)
\]
elements is required for storing the symmetric part. This is the reason why representing tensors not as simple arrays of lists but via decompositions is invaluable for higher-order tensors.

One extensively studied decomposition for tensors is the polyadic decomposition.
\begin{definition}[Polyadic decomposition]
An $n$-dimensional order $d$ tensor $A$ with complex entries is said to have a polyadic decomposition if
\[
A_{i_1,i_2,\ldots,i_d}=\sum_{k=1}^rc_kB_{i_1,k}^{(1)}B_{i_2,k}^{(2)}\cdots B_{i_d,k}^{(d)},\quad i_1,i_2,\ldots,i_d\in\{1,2,\ldots,n\},
\]
where $r\in\mathbb{Z_+}$, $(c_k)_{k=1}^r$ are scalars, and $(B^{(s)})_{s=1}^d$ are $n\times r$ matrices called the \emph{polyadic factors}. 

If the number $r$~is minimal, then the decomposition is called the \emph{canonical polyadic decomposition} or \emph{CP-decomposition}. The smallest possible value for $r$ is called the \emph{canonical rank} or \emph{CP-rank} of tensor $A$ and it is denoted by ${\rm rank}_{\rm CP}A$.
\end{definition}

The polyadic decomposition characterizes the tensor using only $dnr$ parameters, where the complexity of representation transferred into the number of terms $r$.

In the case of symmetric tensors, polyadic decompositions that preserve symmetry are desirable. In this case, the factors of the polyadic decomposition are required to satisfy the additional condition $B^{(1)}=B^{(2)}=\cdots=B^{(d)}$, reducing the the storage complexity to only $nr$, where the smallest possible value for $r$ is called the \emph{symmetric rank} of $A$, ${\rm rank}_{\rm Sym}A$. The symmetric rank is naturally bounded by the canonical rank from below, i.e., ${\rm rank}_{\rm CP}A\leq {\rm rank}_{\rm Sym}A$.

Another type of tensor decomposition considered in this work is the tensor-train decomposition, which was proposed by Oseledets~\cite{Osel11}. The tensor-train decomposition has many desirable qualities from an algorithmic point of view: it can always be formed numerically via an algorithmic sequence of operations, there exists a cross-approximation technique that avoids access to all $n^d$ elements of an $n$-dimensional order $d$ tensor in the construction of the tensor-train decomposition~\cite{OselTyr10}, and if the decomposition possesses a sufficiently low rank, many operations in linear algebra---such as summation, mode contractions against vectors, or dot products with other tensors---can be carried out in polynomial time within the tensor-train framework. The tensor-train decomposition can be seen as a parametric extension of the polyadic decomposition and it is defined as follows. 
\begin{definition}[Tensor-train decomposition]
An $n$-dimensional order $d$ tensor $A$ with complex entries is said to have a tensor-train decomposition or TT-decomposition if
\[
A_{i_1,i_2,\ldots,i_d}=G_1(i_1)G_2(i_2)\cdots G_d(i_d),
\]
where $G_1(i_1)\in\mathbb{C}^{1\times r_1},G_d(i_d)\in\mathbb{C}^{r_{d-1}\times 1}$, and $G_k(i_k)\in\mathbb{C}^{r_{k-1}\times r_k}$ with \emph{compression ranks} $r_1,r_2,\ldots,r_{d-1}\in\mathbb{Z_+}$ for $i_1,i_2,\ldots,i_d\in\{1,2,\ldots,n\}$. The factors $(G_k)_{k=1}^d$ are called the \emph{TT-cores} and they may be regarded as third order tensors by identifying $(G_k)_{j,i,\ell}=G_k(i)_{j,\ell}$.
\end{definition}

The compression ranks are bounded from below by the ranks of the unfolding matrices $A_k$ defined elementwise by
\[
(A_k)_{(i_1,i_2,\ldots,i_k),(i_{k+1},i_{k+2},\ldots,i_d)}=A_{i_1,i_2,\ldots,i_d},
\]
where the multi-indices $(i_1,i_2,\ldots,i_k)\in\{1,2,\ldots,n\}^k$ enumerate the rows and the multi-indices $(i_{k+1},i_{k+2},\ldots,i_d)\in\{1,2,\ldots,n\}^{d-k}$ enumerate the columns of the $n^k\times n^{d-k}$ matrix $A_k$. Moreover, these compression ranks are also achievable.

\begin{theorem}[cf.~{\cite{Osel11}}]\label{osel}
There exists a TT-decomposition of the $n$-dimensional order $d$ tensor $A$ with compression ranks
\[
r_k={\rm rank}\,A_k,\quad k\in\{1,2,\ldots,d-1\}.
\]
\end{theorem}

Theorem~\ref{osel} motivates the definition of the \emph{tensor-train rank} or \emph{TT-rank} of $A$ by setting
\[
{\rm rank}_{\rm TT}A=\max_{1\leq k\leq d-1}{\rm rank}\,A_k.
\]

\subsection{Meet and join tensors}

Let $(P,\preceq,\wedge)$ be a meet semilattice and $f$ a complex-valued function on $P$. The $n$-dimensional order $d$ tensor $(S_d)_f$, where
\[
((S_d)_f)_{i_1,i_2,\ldots,i_d}=f(x_{i_1}\wedge x_{i_2}\wedge\cdots\wedge x_{i_d}),\quad i_1,i_2,\ldots,i_d\in\{1,2,\ldots,n\},
\]
is called the order $d$ \emph{meet tensor} on $S$ with respect to $f$. If $(S,\preceq,\vee)$ is a join semilattice, then the $n$-dimensional order $d$ tensor $[S_d]_f$, where
\[
([S_d]_f)_{i_1,i_2,\ldots,i_d}=f(x_{i_1}\vee x_{i_2}\vee\cdots\vee x_{i_d}),\quad i_1,i_2,\ldots,i_d\in\{1,2,\ldots,n\},
\]
is called the order $d$ \emph{join tensor} on $S$ with respect to $f$. When $d=2$, both of these definitions reduce to the classical and well-studied cases of meet and join matrices.

Let $(P,\preceq,\wedge,\hat 0)$ be a locally finite meet semilattice with the least element $\hat 0$, i.e., $\hat 0\preceq x$ for all $x\in P$, and suppose that $S=\{x_1,x_2,\ldots,x_n\}$ is a finite subset of $P$ such that $x_i\preceq x_j$ only if $i\leq j$. The set $S$ is said to be \emph{meet closed} if $x\wedge y\in S$ for all $x,y\in S$, and $S$ is said to be \emph{join closed} if $x\vee y\in S$ for all $x,y\in S$. Let $f$ be a complex-valued function on $P$.
\[
f_d(\hat 0,x)=f(x),\quad x\in P.
\]

The following $LDL^\textup{T}$ factorization holds for meet matrices on meet closed sets.
\begin{proposition}[cf.~{\cite[Theorem~12]{Ra}}]\label{edet}
Let $S$ be meet closed and define the $n\times n$ matrices $L$ and $D={\rm diag}(d_1,d_2,\ldots,d_n)$ by setting
\begin{align*}
L_{i,j}&=\begin{cases}1,&\text{if }x_j\preceq x_i,\\
0&\text{otherwise},\end{cases}\quad i,j\in\{1,2,\ldots,n\},\notag\\
d_i&=\sum_{\substack{z\preceq x_i\\ z\not\preceq x_j,\ j<i}}(f_d*\mu)(\hat 0,z),\quad i\in\{1,2,\ldots,n\},
\end{align*}
Then $(S_2)_f=LDL^\textup{T}$.
\end{proposition}

This result was generalized by Ilmonen~\cite{ilhyp} for meet tensors on meet closed sets.
\begin{theorem}[cf.~{\cite[Theorem~3.4]{ilhyp}}]
Let $S$ be meet closed and $f$ a complex-valued function on $P$. Then the meet tensor $(S_d)_f$ has a polyadic decomposition given by
\[
((S_d)_f)_{i_1,i_2,\ldots,i_d}=\sum_{k=1}^n d_{k}L_{i_1,k}L_{i_2,k}\cdots L_{i_d,k},\quad i_1,i_2,\ldots,i_d\in\{1,2,\ldots,n\}.
\]
\end{theorem}

When the domain of definition $S$ is not closed, no structure-theoretic results seem to appear in the literature for such classes of lattice-theoretic tensors. In the special case $d=2$ corresponding to matrices, the structure of join matrices is usually related to that of meet matrices on meet closed sets under the assumption that the function $f$ is semimultiplicative~\cite{Is2}. However, we do not make this restriction in the sequel.

\section{Decompositions of join tensors}\label{joindecomp}
A polyadic decomposition of meet tensors was discovered in~\cite{ilhyp}. Tensor decompositions are notoriously difficult to construct and the task of finding the rank of a tensor is generally an NP-hard problem~\cite{Hast90}. Iterative methods such as ALS (a detailed description of which can be found, e.g., in the survey~\cite{KolBa09}) can sometimes be used to derive numerical decompositions, but their convergence is not guaranteed---not to mention that numerical decompositions are usually of little theoretical interest.
\subsection{Construction of a polyadic decomposition of join tensors}\label{constrcp}
We present in the following one approach to construct a polyadic decomposition for join tensors. The decomposition can be expressed using lattice-theoretical quantities and as such may be theoretically meaningful.

Let $(P,\preceq,\vee)$ be a locally finite join semilattice and $S=\{x_1,x_2,\ldots,x_n\}$ a finite subset of $P$. In the sequel, we denote 
\[
S^{\vee k}=\{x_{i_1}\vee x_{i_2}\vee\cdots\vee x_{i_k}\mid 1\leq i_1\leq i_2\leq\cdots\leq i_k\leq n\}\quad\text{for }k\in\mathbb{Z_+}.
\]
The approach taken in this paper can be seen as a generalization of the methods developed in the analysis of join matrices defined over two posets~\cite{MaHa13}, where the set $S^{\vee 2}$ plays an important role.
\begin{theorem}\label{cpdecomp}
Let $(P,\preceq,\vee)$ be a locally finite join semilattice, $S=\{x_1,x_2,\ldots,x_n\}$ a finite subset of $P$ ordered such that $x_i\preceq x_j$ only if $i\leq j$, and $f$ a complex-valued function on $P$. Let $S^{\vee d}=\{y_1,y_2,\ldots,y_r\}$ be ordered such that $y_i\preceq y_j$ only if $i\leq j$. Then $[S_d]_f$ has a $\# S^{\vee d}$-term polyadic decomposition given by
\[
([S_d]_f)_{i_1,i_2,\ldots,i_d}=\sum_{k=1}^rc_kE_{i_1,k}E_{i_2,k}\cdots E_{i_d,k},
\]
where
\begin{align}
c_k=\sum_{s:\,y_k\preceq y_s}f(y_s)\mu_{S^{\vee d}}(y_k,y_s),\quad k\in\{1,2,\ldots,r\}.\label{cpcoeffs}
\end{align}
Here $\mu_{S^{\vee d}}$ denotes the M\"{o}bius function of the poset $(S^{\vee d},\preceq)$ and $E$ is the $n\times r$ matrix defined elementwise by setting
\[
E_{i,j}=\begin{cases}1,&\text{if }x_i\preceq y_j,\\ 0&\text{otherwise}\end{cases}
\]
for $i\in\{1,2,\ldots,n\}$ and $j\in\{1,2,\ldots,r\}$.
\end{theorem}
\proof Let the scalars
\[
c_k=\sum_{s:\,y_k\preceq y_s}f(y_s)\mu_{S^{\vee d}}(y_k,y_s),\quad k\in\{1,2,\ldots,r\}.
\]
By the dual formulation of the M\"{o}bius inversion formula~\cite[Proposition~3.7.2]{St} this is equivalent to
\[
f(y_s)=\sum_{k:\,y_s\preceq y_k}c_k=\sum_{k=1}^rc_k\zeta(y_s,y_k),
\]
where the latter equality is valid since $\zeta(y_s,y_k)=0$ for all $y_k\prec y_s$ thereby removing the dependence on $y_s$ from the summation variable.

Since the set $S^{\vee d}$ contains, by definition, all of the unique elements of the tensor $[S_d]_f$, then for each tuplet $(x_{i_1},\ldots,x_{i_d})\in S^d$ there exists a representative $y_s\in S^{\vee d}$ such that $y_s=x_{i_1}\vee\cdots\vee x_{i_d}$. Hence this equality can be substituted into the formula above to obtain
\begin{align*}
f(x_{i_1}\vee\cdots\vee x_{i_d})=\sum_{k=1}^rc_k\zeta(x_{i_1}\vee\cdots\vee x_{i_d},y_k)
\end{align*}
and since the join operation is the least upper bound of two elements, it immediately follows from recursive application of the universal property $x,y\preceq z\Leftrightarrow x\vee y\preceq z$ for any $x,y,z\in P$ that
\begin{align*}
f(x_{i_1}\vee\cdots\vee x_{i_d})&=\sum_{k=1}^rc_k\zeta(x_{i_1}\vee\cdots\vee x_{i_d},y_k)\\
&=\sum_{k=1}^rc_k\zeta(x_{i_1},y_k)\cdots\zeta(x_{i_d},y_k)\\
&=\sum_{k=1}^rc_kE_{i_1,k}\cdots E_{i_d,k},
\end{align*} 
proving the assertion.\endproof
\begin{remark}
In practice, the computation of the coefficient vector $c=(c_1,\ldots,c_r)^\textup{T}$ can be carried out by solving the matrix equation
\begin{align}
\mathcal{E}c=y,\label{matrixequ}
\end{align}
where $y=(f(y_1),\ldots,f(y_r))^\textup{T}$ and $\mathcal{E}$ is the $r\times r$ matrix defined elementwise by setting
\[
\mathcal{E}_{i,j}=\begin{cases}1,&\text{if }y_i\preceq y_j,\\ 0&\text{otherwise}\end{cases}
\]
for $i,j\in\{1,\ldots,n\}$. The matrix $\mathcal{E}$ is nonsingular since it is an upper triangular matrix and its diagonal contains no zeros. In fact, its inverse can be given elementwise by $\mathcal{E}_{i,j}^{-1}=\mu_{S^{\vee d}}(y_i,y_j)$, which makes the relationship between the matrix equation~\eqref{matrixequ} and the expression of the coefficients in~\eqref{cpcoeffs} obvious. In other words, solving the matrix equation~\eqref{matrixequ} is equivalent to M\"{o}bius inversion.
\end{remark}

Determining the elements of the set $S^{\vee d}$ may require up to $\binom{n+d-1}{d}=\mathcal{O}(n^d/d!)$
evaluations (see the discussion below Remark~\ref{symremark}). The polyadic decomposition of Theorem~\ref{cpdecomp} may be interpreted as a compressed representation of the symmetric tensor $[S_d]_f$. That is, if the unique elements of $[S_d]_f$ have already been computed, then the full tensor may be represented economically via the polyadic decomposition, requiring only the storage of the $r$ coefficients $(c_k)_{k=1}^r$ and the $n\times r$ matrix $E$, the structure (and potential sparsity) of which depends on the semilattice.

The worst-case scenario with the largest amount of stored elements is realized if the join operator satisfies
\begin{align}
x_{p_1}\vee\cdots\vee x_{p_d}=x_{q_1}\vee\cdots\vee x_{q_d}\quad\text{if and only if}\quad \{p_1,\ldots,p_d\}=\{q_1,\ldots,q_d\}\label{worst}
\end{align}
for all elements of $S=\{x_1,x_2,\ldots,x_n\}\subset P$ when $f\!:P\to\mathbb{C}$ is a one-to-one mapping. The cardinality of the set $S^{\vee k}$ for $k\in\mathbb{Z}_+$ (sic) can be bounded by
\[
\# S^{\vee k}\leq \# S^{\vee n}\leq\sum_{i=1}^{\# S}\binom{\# S}{i}=2^{\# S}-1,
\]
where the upper bound is achieved only under the condition~\eqref{worst}. 

In the average setting---such as in the case of the divisor lattice $(\mathbb{Z_+},|,{\rm gcd},{\rm lcm})$---the cardinality $\# S^{\vee d}$ may be substantially less than the worst-case cardinality. In the following example, we study the profile of the polyadic factors of a class of LCM tensors for increasing dimension $n$.

\begin{example}
We examine the polyadic factors of the LCM tensor $[(S^{(n)})_8]$ defined on the locally finite divisor lattice $(\mathbb{Z_+},|,{\rm gcd},{\rm lcm})$ with respect to $S^{(n)}=\{1,2,\ldots,n\}$ by setting
\[
([(S^{(n)})_8])_{i_1,i_2,\ldots,i_8}={\rm lcm}(i_1,i_2,\ldots,i_8),\quad i_1,i_2,\ldots,i_8\in\{1,\ldots,n\},
\]
where we let $n\in\{1,2,\ldots,12\}$ and denote the associated polyadic factors given by Theorem~\ref{cpdecomp} with $E^{(n)}$. For ease of presentation, the columns of the matrix $E^{(12)}$ are ordered in such a way that $E^{(n)}$ is always a leading submatrix of $E^{(n+1)}$ for $n\in\{1,\ldots,11\}$. The profiles of the polyadic factors are displayed in Figure~\ref{mtxplot}, where each $E^{(n)}$ is sectioned off for $n\in\{1,\ldots,12\}$ counting from top left to bottom right. Since each polyadic factor is a Boolean matrix, the matrix plot is presented as a grid where gray coloring denotes elements equal to unity and white coloring denotes elements equal to zero.

\begin{figure}[!h]
\centering
\includegraphics[width=14cm]{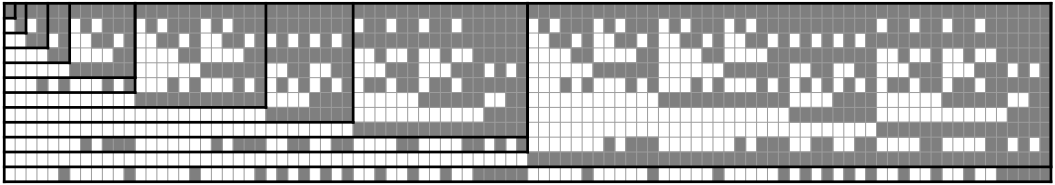}
\caption{The matrix plot of the polyadic factors $E^{(n)}$ of the LCM tensor $[(S^{(n)})_8]$, where each $E^{(n)}$ is sectioned off for $n\in\{1,\ldots,12\}$. The gray coloring denotes matrix elements equal to unity and white coloring denotes elements equal to zero.}\label{mtxplot}
\end{figure}

Additional numerical experiments on the number of nonzero elements in the polyadic factors given by Theorem~\ref{cpdecomp} are carried out in Subsection~\ref{numex}.

\end{example}

We also mention that in the special case when the set $S=\{x_1,x_2,\ldots,x_n\}\subset P$ is join closed, the number of terms of the decomposition of Theorem~\ref{cpdecomp} is at most $n$. This case is realized, e.g., in the case of the locally finite lattice $(\mathbb{Z_+},\leq,\min,\max)$ with respect to $S=\{1,2,\ldots,n\}$, where the associated MAX tensor is defined by setting
\[
([S_d])_{i_1,i_2,\ldots,i_d}=\max\{i_1,i_2,\ldots,i_d\},\quad i_1,i_2,\ldots,i_d\in\{1,2,\ldots,n\}.
\]

Polyadic decompositions may be useful in the assessment of theoretical properties of join tensors, such as in the analysis of tensor eigenvalues. However, this field of study is still under development and requires further research. Moreover, the difficulty in the numerical construction of the polyadic decomposition of Theorem~\ref{cpdecomp} limits its practical usage. We next move onto the tensor-train decomposition of join tensors, which provides one way to overcome the computational load associated with the computation of the set $S^{\vee d}$.

\subsection{Construction of a tensor-train decomposition of join tensors}
The motivation behind the tensor-train decomposition stems from the fact that it is easy to construct numerically in an algorithmic sequence of operations for practically any tensor~\cite{Osel11}. In the case of join tensors, however, it turns out that an explicit TT-decomposition can be derived.

The explicit TT-decomposition we derive utilizes the inherent symmetry of join tensors. Additionally, the TT-cores of this representation are remarkably sparse. Neither of these features arise by constructing the TT-decomposition numerically via the TT-SVD or TT-Cross algorithms (see~\cite{OselTyr10}), which makes the explicit TT-decomposition of higher-order join tensors preferable from the viewpoint of numerical analysis. In addition, the TT-rank is easy to analyze since it can be reduced to inspecting the ranks of unfolding matrices.

\begin{figure}[!h]
\begin{center}
\includegraphics{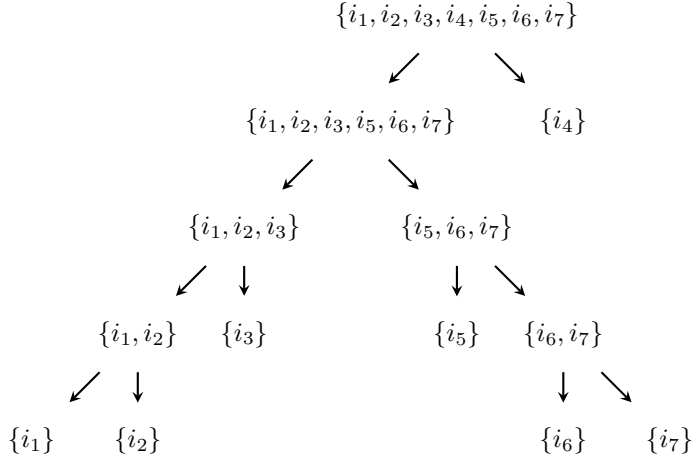}
\caption{An illustration of the partition tree for $[S_7]_f$. The center mode is separated at the first step. From the second step onward, the remaining modes can be split using a linear partition tree.}\label{dimtree}
\end{center}
\end{figure}

The ultimate goal of the tensor-train decomposition is to separate the indices of the tensor. By carrying out this partitioning carefully, it is possible to reduce the computational load associated with the construction of the decomposition. The partition tree associated with $d=7$ is illustrated as an example in Figure~\ref{dimtree}. For general $d$, we separate the $(\lfloor d/2\rfloor+1)^\text{th}$ mode at the first step and proceed to separate the remaining modes via a linear partition tree. The conditions for the optimality of this construction are discussed in Section~\ref{ttrank}.

The structure theorem is presented in the following.
\begin{theorem}\label{ttdecomp}
Let $(P,\preceq,\vee)$ be a locally finite join semilattice, $S=\{x_1,x_2,\ldots,x_n\}$ a finite subset of $P$ ordered such that $x_i\preceq x_j$ only if $i\leq j$, and $f$ a complex-valued function on $P$. Then there exists a tensor-train decomposition with compression ranks at most $\# S^{\vee\lfloor d/2 \rfloor}$ of the $n$-dimensional order $d$ join tensor $[S_d]_f$ given by
\[
([S_d]_f)_{i_1,i_2,\ldots,i_d}=G_1(i_1)G_2(i_2)\cdots G_d(i_d),\quad i_1,i_2,\ldots,i_d\in\{1,2,\ldots,n\},
\]
where the TT-cores are defined by setting
\begin{align*}
&G_1(i_1)_{1,\alpha_1}=\delta(x_{i_1},\alpha_1),\\
&G_k(i_k)_{\alpha_{k-1},\alpha_k}=\delta(\alpha_k,x_{i_k}\vee \alpha_{k-1}),\quad k\in\left\{2,3,\ldots,\left\lfloor\frac{d}{2}\right\rfloor\right\},\\
&G_{\lfloor d/2\rfloor+1}(i_{\lfloor d/2\rfloor+1})_{\alpha_{\lfloor d/2\rfloor},\alpha_{d-\lfloor d/2\rfloor-1}}=f(\alpha_{\lfloor d/2\rfloor}\vee x_{i_{\lfloor d/2\rfloor+1}}\vee\alpha_{d-\lfloor d/2\rfloor-1}),\\
&G_k(i_k)_{\alpha_{d-k+1},\alpha_{d-k}}=G_{d-k+1}(i_k)_{\alpha_{d-k},\alpha_{d-k+1}},\quad k\in\left\{\left\lfloor\frac{d}{2}\right\rfloor+2,\ldots,d-1\right\},\\
&G_d(i_d)_{\alpha_1,1}=G_1(i_d)_{1,\alpha_1},
\end{align*}
where $\alpha_k\in S^{\vee k}$ for $k\in\{1,2,\ldots,\lfloor d/2\rfloor\}$.
\end{theorem}
\proof Let us first consider join tensors of even order $d=2q$, $q\geq 1$. Since $([S_{d}]_f)_{i_1,\ldots,i_{d}}=f(x_{i_1}\vee\cdots\vee x_{i_{d}})$, it is possible to expand this as
\[
([S_{d}]_f)_{i_1,\ldots,i_{d}}=\sum_{\alpha_{q},\alpha_{q-1}} \delta(\alpha_{q},x_{i_1}\vee\cdots\vee x_{i_{q}}) f(\alpha_{q}\vee x_{i_{q+1}}\vee \alpha_{q-1}) \delta(\alpha_{q-1},x_{i_{q+2}}\vee\cdots\vee x_{i_{2q}})
\]
where the sum is taken over $\alpha_{q}\in S^{\vee q}$ and $\alpha_{q-1}\in S^{\vee (q-1)}$. The claim follows by observing that
\[
\delta(\alpha_{q},x_{i_1}\vee\cdots\vee x_{i_{q}})=\sum_{\alpha_1,\ldots,\alpha_{q-1}} \delta(\alpha_1,x_{i_1})\delta(\alpha_2,\alpha_1\vee x_{i_2})\cdots \delta(\alpha_{q},\alpha_{q-1}\vee x_{i_{q}})
\]
is in the desired tensor-train format when the sum is taken over each $\alpha_k\in S^{\vee k}$ for $k\in\{1,2,\ldots,q-1\}$ and, likewise,
\begin{align*}
&\delta(\alpha_{q-1},x_{i_{q+2}}\vee\cdots\vee x_{i_{2q}})\\
&=\sum_{\alpha_{q-2},\ldots,\alpha_{1}} \delta(\alpha_{q-1},\alpha_{q-2}\vee x_{i_{q+2}})\cdots \delta(\alpha_{2},\alpha_{1}\vee x_{i_{2q-1}})\delta(\alpha_{1},x_{i_{2q}}),
\end{align*}
where the sum is taken over each $\alpha_{k}\in S^{\vee k}$ for $k\in\{1,2,\ldots,q-2\}$.

Next, we consider join tensors of odd order $d=2q+1$, $q\geq 1$. This time we write
\begin{align*}
&([S_{d}]_f)_{i_1,\ldots,i_{d}}\\
&=\sum_{\alpha_q,\alpha_{q}'} \delta(\alpha_q,x_{i_1}\vee\cdots\vee x_{i_q})f(\alpha_q\vee x_{i_{q+1}}\vee \alpha_{q}')\delta(\alpha_{q}',x_{i_{q+2}}\vee\cdots\vee x_{i_{2q+1}}),
\end{align*}
where the sum is taken over $\alpha_q,\alpha_{q}'\in S^{\vee q}$. The rest of the proof proceeds as in the even order case.\endproof
\begin{remark}
The elements of $S^{\vee k}$ have been used to enumerate the rows and columns of the TT-factors in Theorem~\ref{ttdecomp} for ease of presentation. In practice, it is trivial to construct a renumbering $\tau_k\!:S^{\vee k}\ni\alpha_k\mapsto i\in\{1,\ldots,\# S^{\vee k}\}$ for $k\in\{1,\ldots,\lfloor d/2\rfloor\}$ such that $G_k(i_k)_{\alpha_{k-1},\alpha_k}=G_k(i_k)_{\tau_{k-1}(\alpha_{k-1}),\tau_{k-1}(\alpha_k)}$.
\end{remark}
\begin{remark}
In the locally finite divisor lattice $(\mathbb{Z_+},|,{\rm gcd},{\rm lcm})$ with $S=\{1,2,\ldots,n\}$, there exists a natural  interpretation for the sets $S^{\vee k}$. Namely, 
\begin{align*}
S^{\vee k}&=\{{\rm lcm}(i_1,\ldots,i_k)\mid 1\leq i_1\leq\cdots\leq i_k\leq n\}\\
&=\{i_1\cdots i_k\mid 1\leq i_j\leq n\text{ and }{\rm gcd}(i_p,i_q)=1,\ p\neq q\},
\end{align*}
which follows from the prime number decomposition of ${\rm lcm}(i_1,\ldots,i_k)$ and the fact that ${\rm lcm}(i_1,\ldots,i_k)=i_1\cdots i_k$ if and only if ${\rm gcd}(i_p,i_q)=1$ for all $p\neq q$. This special case has been discussed in a separate work~\cite{hik16}.
\end{remark}

The utility of the tensor-train decomposition stems from the fact that only the first $\lfloor d/2\rfloor+1$ TT-cores need to be stored explicitly. The remaining TT-cores are nothing more than transposes of the first cores in reverse order. In practice, the storage cost of the tensor-train decomposition is dominated by the $(\lfloor d/2\rfloor+1)^\text{th}$ TT-core since all other TT-cores are sparse, Boolean tensors.

Unlike the polyadic decomposition, the TT-decomposition can be formed efficiently for higher-order join tensors as it requires only the construction of the sets $S^{\vee k}$ for $k\in\{2,\ldots,\lfloor d/2\rfloor\}$, which is significantly less computationally taxing compared to the construction of the set $S^{\vee d}$ in the case of the polyadic decomposition.

\section{On the rank of join tensors}\label{ttrank}
The structure of join tensors allows for the analysis of the TT-rank of join tensors, which is carried out in Subsection~\ref{ttranksub}. This, in turn, can be used to derive a lower bound for the CP-rank via a straightforward embedding argument provided in Subsection~\ref{cpranksub}. It is also possible to show that for certain join tensors satisfying $n\leq\lfloor d/2\rfloor$, the explicit polyadic decomposition presented in Theorem~\ref{cpdecomp} is canonical and, under the same condition, the tensor-train decomposition of Theorem~\ref{ttdecomp} is a TT-rank decomposition.

\subsection{On the TT-rank of join tensors}\label{ttranksub}
We first investigate the bounds on the TT-rank of join tensors.
\begin{theorem}\label{ranknd}
Let $(P,\preceq,\vee)$ be a locally finite join semilattice, $S=\{x_1,x_2,\ldots,x_n\}$ a finite subset of $P$ such that $x_i\preceq x_j$ only if $i\leq j$, and $f\!:P\to\mathbb{C}$ a function on $P$. Let $S^{\vee d}=\{y_1,y_2,\ldots,y_{r}\}$ be ordered $y_i\preceq y_j$ only if $i\leq j$ and assume that $f$ satisfies
\begin{align}
c_k=\sum_{s:\,y_k\preceq y_s}f(y_s)\mu_{S^{\vee d}}(y_k,y_s)\neq 0\quad\text{for all }k\in\{1,2,\ldots,r\}.\label{coeffassumption}
\end{align}
Then the TT-rank of the $n$-dimensional order $d$ join tensor $[S_d]_f$ is bounded by
\[
2\# S^{\vee \lfloor d/2\rfloor}-\# S^{\vee d}\leq {\rm rank}_{\rm TT}[S_d]_f\leq\#S^{\vee \lfloor d/2\rfloor}.
\]
\end{theorem}
\proof
The upper bound ${\rm rank}_{\rm TT}[S_d]_f\leq \#S^{\vee\lfloor d/2\rfloor}$ follows immediately from the construction of Theorem~\ref{ttdecomp}. Finding a lower bound for the TT-rank is the nontrivial part of the argument.

The $k^\text{th}$ unfolding matrix of the join tensor $A=[S_d]_f$ is given by
\[
(A_k)_{(i_1,\ldots,i_k),(i_{k+1},\ldots,i_d)}=f((x_{i_1}\vee\cdots\vee x_{i_k})\vee(x_{i_{k+1}}\vee\cdots\vee x_{i_d})),
\]
and to establish the lower bound, we select a $\#S^{\vee \lfloor d/2\rfloor}\times\#S^{\vee \lfloor d/2\rfloor}$ submatrix $B$ of $A_{\lfloor d/2\rfloor}$ defined by setting
\[
B_{\alpha,\beta}=f(\alpha\vee\beta),
\]
where $\alpha,\beta\in S^{\vee \lfloor d/2\rfloor}$ are used to enumerate the rows and columns of $B$, respectively.\footnotemark\footnotetext{For $d$ even, this submatrix exists trivially. For $d$ odd, one uses the fact that $S^{\vee \lfloor d/2\rfloor}\subset S^{\vee\lceil d/2\rceil}$ since any element $x_{i_1}\vee\cdots \vee x_{i_d}\in S^{\vee\lfloor d/2\rfloor}$ is also an element $x_{i_1}\vee\cdots\vee x_{i_{\lfloor d/2\rfloor}}\vee x_{i_k}\in S^{\vee \lceil d/2\rceil}$ for any $1\leq k\leq \lfloor d/2\rfloor$.} This submatrix can be interpreted as a join matrix over the poset $(S^{\vee \lfloor d/2\rfloor},\preceq)$ and as such it has a decomposition given by $B=\Gamma\Lambda \Gamma^\textup{T}$ (see~\cite[Theorem~3.1]{MaHa13}), where $\Lambda={\rm diag}(c_1,\ldots,c_r)$ and $\Gamma$ is a $\# S^{\vee \lfloor d/2\rfloor}\times \#S^{\vee d}$ matrix defined elementwise by
\[
\Gamma_{i,j}=\begin{cases}1,&\text{if }z_i\preceq y_j,\\ 0&\text{otherwise},\end{cases}
\]
where $S^{\vee \lfloor d/2\rfloor}=\{z_1,\ldots,z_{\# S^{\vee\lfloor d/2\rfloor}}\}$ is ordered such that $z_i\preceq z_j$ only if $i\leq j$.

The matrix factor $\Gamma$ has full rank since $S^{\vee \lfloor d/2\rfloor}\subset S^{\vee d}$. On the other hand, the assumption~\eqref{coeffassumption} implies that the diagonal matrix $\Lambda$ also has full rank. Since $B$ is a submatrix of $A_{\lfloor d/2\rfloor}$, it now holds that ${\rm rank}\,B\leq{\rm rank}\,A_{\lfloor d/2\rfloor}$ and applying Sylvester's rank inequality\footnotemark\footnotetext{Sylvester's rank inequality: ${\rm rank}(AB)\geq {\rm rank}\,A+{\rm rank}\,B-n$ for any $k\times n$ matrix $A$ and $n\times m$ matrix $B$.} to the matrix product $\Gamma\Lambda\Gamma^\textup{T}$ yields
\begin{align*}
{\rm rank}\,A_{\lfloor d/2\rfloor}&\geq{\rm rank}\,B\\
&={\rm rank}(\Gamma\Lambda\Gamma^{\textup{T}})\\
&\geq {\rm rank}\,\Gamma+{\rm rank}(\Lambda\Gamma^{\textup{T}})-\# S^{\vee d}\\
&\geq {\rm rank}\,\Gamma+{\rm rank}\,\Lambda+{\rm rank}\,\Gamma-2\,\#S^{\vee d}\\
&=2\# S^{\vee \lfloor d/2\rfloor}-\# S^{\vee d}.
\end{align*}
In particular, it follows that
\[
{\rm rank}_{\rm TT}A=\max_{1\leq k\leq d-1}{\rm rank}\,A_k\geq {\rm rank}\,A_{\lfloor d/2\rfloor}\geq 2\# S^{\vee \lfloor d/2\rfloor}-\# S^{\vee d}.
\]
This proves the assertion.\endproof

As a corollary, we obtain the following result.
\begin{corollary}\label{ttcoro}
If in addition to the assumptions of Theorem~\ref{ranknd} it holds that $n\leq\lfloor d/2\rfloor$, then
\[
{\rm rank}_{\rm TT}[S_d]_f=\# S^{\vee\lfloor d/2\rfloor}.
\]
\end{corollary}
\proof It is easy to see by the pigeon hole principle that $S^{\vee n}=S^{\vee k}$ for all $k\geq n$. If $n\leq\lfloor d/2\rfloor$, then it immediately follows that $S^{\vee \lfloor d/2\rfloor}=S^{\vee d}$.
\endproof

\subsection{On the CP-rank of join tensors}\label{cpranksub}
A trivial but nonetheless important consequence of Theorem~\ref{ranknd} is the fact that it also provides a lower bound for the CP-rank of tensor $[S_d]_f$ as well. If $A=[S_d]_f$ has a CP-decomposition of the form
\[
A_{i_1,i_2,\ldots,i_d}=\sum_{\alpha=1}^rc_\alpha B_{i_1,\alpha}^{(1)}B_{i_2,\alpha}^{(2)}\cdots B_{i_d,\alpha}^{(d)},
\]
where $B^{(k)}\in\mathbb{C}^{n\times r}$, $c_\alpha\in\mathbb{C}$, and $r={\rm rank}_{\rm CP}A$, then the canonical polyadic form can be embedded into the tensor-train format by
\begin{align*}
&\sum_{\alpha=1}^rc_\alpha B_{i_1,\alpha}^{(1)}B_{i_2,\alpha}^{(2)}\cdots B_{i_d,\alpha}^{(d)}\\
&=\sum_{\alpha_1,\ldots,\alpha_{d-1}=1}^rc_{\alpha_1} B_{i_1,\alpha_1}^{(1)}\cdot \delta(\alpha_1,\alpha_2)B_{i_2,\alpha_2}^{(2)}\cdots \delta(\alpha_{d-2},\alpha_{d-1})B_{i_{d-1},\alpha_{d-1}}^{(d-1)}\cdot B_{i_d,\alpha_{d-1}}^{(d)},
\end{align*}
where the latter term is in tensor-train formalism. This implies the following.
\begin{corollary}\label{cprank} Under the assumptions of Theorem~\ref{ranknd}, it holds that
\[
2\# S^{\vee\lfloor d/2\rfloor}-\# S^{\vee d}\leq {\rm rank}_{\rm CP}[S_d]_f\leq{\rm rank}_{\rm Sym}[S_d]_f\leq\#S^{\vee d}.
\]
If $n\leq\lfloor d/2\rfloor$, then
\[
{\rm rank}_{\rm CP}[S_d]_f={\rm rank}_{\rm Sym}[S_d]_f=\#S^{\vee \lfloor d/2\rfloor}.
\]
\end{corollary}
\proof The lower bound follows from the embedding of the CP-decomposition into TT-format. Conversely, the upper bound holds because of the explicit factorization given by Theorem~\ref{cpdecomp}.

The second part follows analogously to the proof of Corollary~\ref{ttcoro}.\endproof

\section{Applications to tensor eigenvalue problems}\label{ttappl}
The decompositions presented in Section~\ref{joindecomp} provide low parametric representations of join tensors and are useful in applications involving high order and dimensionality. In this section, we investigate the storage complexity of these representations. Since the tensor-train and polyadic decompositions allow for fast linear algebraic operations with join tensors, the obtained decompositions can be used for the solution of dominant eigenvalues of join tensors. To support the numerical computations, an elementary upper bound on the eigenvalues of join tensors is presented. This upper bound has appeared recently in~\cite[Theorem~4.1]{ilhyp} and we state it in the case when the assumption of closeness of the domain of definition is removed. The sharpness of this upper bound is assessed for higher-order LCM tensors, which are a special case of join tensors. 
\subsection{An upper bound on the eigenvalues of join tensors}
The tensor eigenvalue problem was introduced by Qi~\cite{Qi} and Lim~\cite{Lim} and it has been considered for meet tensors on meet closed sets by Ilmonen~\cite{ilhyp}. The number $\lambda\in\mathbb{C}$ is an eigenvalue if it satisfies
\begin{align}
Ax^{d-1}=\lambda x^{[d-1]}\label{tevp}
\end{align}
for some nonzero eigenvector $x=(x_1,\ldots,x_n)^\textup{T}$.

In the following, the extended Kronecker symbol defined by
\[
\delta(i_1,i_2,\ldots,i_d)=\begin{cases}1,&\text{if }i_1=i_2=\cdots=i_d,\\ 0&\text{otherwise}\end{cases}
\]
is used for notational convenience.
\begin{proposition}[Generalization of {\cite[Theorem~4.1]{ilhyp}}]\label{evpbound}
Let $(P,\preceq,\vee)$ be a locally finite join semilattice, $S=\{x_1,x_2,\ldots,x_n\}$ a finite subset of $P$ ordered $x_i\preceq x_j$ only if $i\leq j$, and let $f$ be any real-valued function on $P$. Then every eigenvalue of the $n$-dimensional order $d$ join tensor $[S_d]_f$ lies in the region
\[
\bigcup_{i=1}^n\{z\in\mathbb{C}: |z-f(x_i)|\leq (n^{d-1}-1)c_i\},
\]
where we define
\begin{align}
c_i=\max_{\substack{i_2,\ldots,i_d\in\{1,2,\ldots,n\}\\ \delta(i,i_2,\ldots,i_d)=0}}|f(x_i\vee x_{i_2}\vee\cdots\vee x_{i_d})|\quad\text{for }i\in\{1,2,\ldots,n\}.\label{ci}
\end{align}
\end{proposition}
\proof
The proof follows the same steps as in~\cite[Theorem~4.1]{ilhyp}, where an upper bound is obtained by utilizing the Gerschgorin theorem for tensors~\cite{Qi}. One needs to use the uniform upper bound~\eqref{ci} when the set $S$~is not closed under the join (resp. meet) operation.\endproof
\begin{remark}
Naturally, the upper bound is affected by the nature of the lattice and associated incidence function. The following corollaries follow immediately from Proposition~\ref{evpbound}.
\begin{itemize}
\item[(i)] If $S$ is join closed, then $c_k\leq\max_{1\leq i\leq n}|f(x_i)|$ for all $k\in\{1,2,\ldots,n\}$ and the result is equivalent to~\cite[Theorem~4.1]{ilhyp}.
\item[(ii)] If $f$ is real-valued and absolutely decreasing in $P$, i.e., $x\preceq y$ implies that $|f(y)|\leq |f(x)|$ for all $x\in P$, then $c_k=f(x_k)$ for all $k\in\{1,2,\ldots,n\}$.
\end{itemize}
\end{remark}
\subsection{Numerical experiments}\label{numex}
In our numerical experiments, we consider a class of LCM tensors on the locally finite divisor lattice $(\mathbb{Z_+},|,{\rm gcd},{\rm lcm})$ with respect to $S=\{1,2,\ldots,n\}$ defined elementwise by setting
\[
([S_d])_{i_1,i_2,\ldots,i_d}={\rm lcm}(i_1,i_2,\ldots,i_d),\quad i_1,i_2,\ldots,i_d\in\{1,2,\ldots,n\}.
\]
The numerical experiments have been divided into two parts. In Subsection~\ref{sparsity} we investigate the storage complexity of the explicit polyadic and TT-representations of $[S_d]$ and the number of nonzero elements of these representations is compared to the number of elements of the associated symmetric part ${\rm Sym}\,[S_d]$. In Subsection~\ref{eigenvaluecomputations}, we compute the dominant eigenvalues of $[S_d]$ for an ensemble of $d$ and $n$ and compare the obtained values to the theoretical upper bound given by Proposition~\ref{evpbound}.
\subsubsection{Low parametric representation of LCM tensors in tensor-train and polyadic formats}\label{sparsity}
\begin{figure}[!h]
\centering
\includegraphics{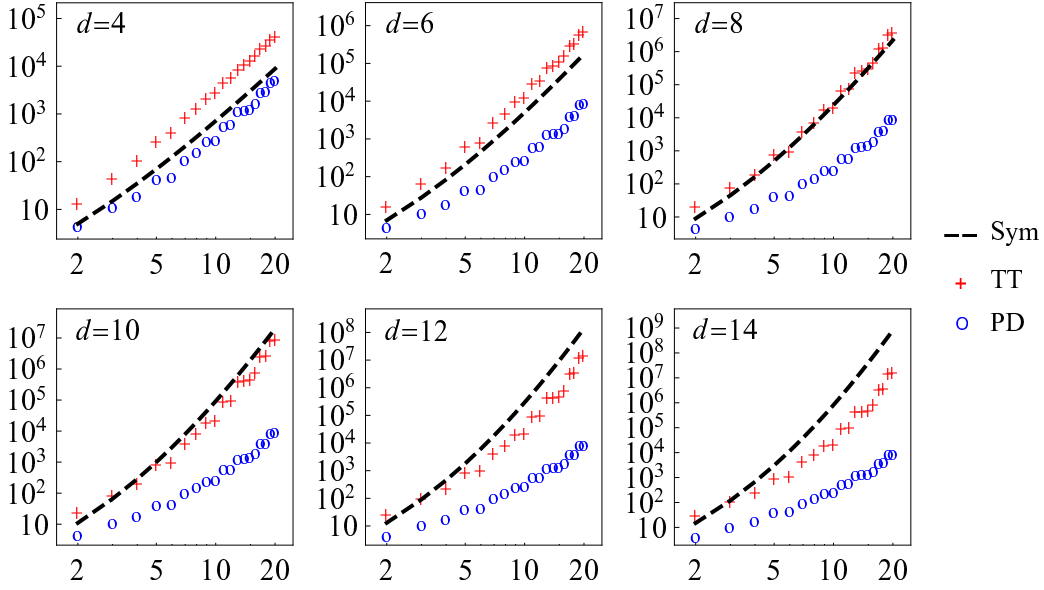}
\caption{The number of nonzero elements that need to be stored in the symmetric part, TT-decomposition, and polyadic decomposition of the LCM tensor $[S_d]$, $S=\{1,\ldots,n\}$, with increasing dimension $2\leq n\leq 20$ and fixed order $d\in\{4,6,8,10,12,14\}$.}\label{sparsityfig}
\end{figure}

The TT-decomposition of $[S_d]$ was computed using Theorem~\ref{ttdecomp}. As such, the storage complexity of the TT-decomposition is characterized by the number of nonzero elements in the TT-cores $(G_k)_{k=1}^{\lfloor d/2\rfloor+1}$. On the other hand, the polyadic decomposition was computed using Theorem~\ref{cpdecomp} and it is parametrized by the number of nonzero elements in the factor matrix $E$ and the coefficient vector of length $\#S^{\vee d}$. We compare the number of these parameters to the number of elements of the symmetric part ${\rm Sym}\,[S_d]$, which is equal to $\binom{d+n-1}{d}$. The results are displayed in Figure~\ref{sparsityfig}.

The polyadic decomposition is superior in its sparsity compared to either the TT-decomposition or the symmetric part regardless of order or dimensionality. However, the difficulty in forming the polyadic decomposition makes its construction intractable for orders $d>15$ and dimensionality $n>20$ in practice. The TT-decomposition can be formed consistently for orders roughly twice that of the comparable polyadic decomposition. The TT-decomposition is parametrized by a slightly greater number of elements than the symmetric part for $d<10$. For $d\geq 10$, however, the storage complexity of the TT-decomposition overtakes that of the symmetric part.
\subsubsection{Solution of dominant eigenvalues of LCM tensors}\label{eigenvaluecomputations}
The dominant eigenvalues of a class of symmetric tensors can be solved by using the following algorithm.
\begin{alg}[Power method, cf.~{\cite{Ng09}}]\label{power}Let $A$ be a symmetric, $n$-dimensional even order $d$ tensor with positive elements.
\begin{itemize}[label={}]
\item Set initial guess $x^{(0)}\in\mathbb{R}_+^n$.\\
For $k=1,2,\ldots$ do
\begin{itemize}[label={}]
\item Compute $y^{(k)}=A(x^{(k-1)})^{d-1}$.
\item Set $x^{(k)}=\displaystyle\frac{(y^{(k)})^{\left[\frac{1}{d-1}\right]}}{\big\|(y^{(k)})^{\left[\frac{1}{d-1}\right]}\big\|}$.
\item Compute $\underline{\lambda}_k=\displaystyle\min_{x_i^{(k)}>0}\displaystyle\frac{(A(x^{(k)})^{d-1})_i}{(x_i^{(k)})^{d-1}}$ and $\overline{\lambda}_k=\displaystyle\max_{x_i^{(k)}>0}\displaystyle\frac{(A(x^{(k)})^{d-1})_i}{(x_i^{(k)})^{d-1}}$.
\end{itemize}
until $\overline{\lambda_k}-\underline{\lambda_k}<\varepsilon$.\\
The dominant eigenvalue $\lambda$ of $A$ satisfies
\[
\underline{\lambda_1}\leq\cdots\leq\underline{\lambda_k}\leq\cdots\leq\lambda\leq\cdots\leq \overline{\lambda_k}\leq\cdots\leq\overline{\lambda_1}.
\]
\end{itemize}
\end{alg}

We make the following observations that reduce the computational complexity of Algorithm~\ref{power}. If the tensor $A$ is characterized by the TT-cores $(G_k)_{k=1}^d$, then the contractions $Ax^{d-1}$ and $Ax^d$ can be carried out efficiently using the formulae
\begin{align}
&Ax^{d-1}=G_1(G_2\times_2 x)(G_3\times_2 x)\cdots (G_d\times_2 x),\label{contr1}\\
&Ax^d=(G_1\times_2 x)(G_2\times_2 x)\cdots (G_d\times_2 x),\label{contr2}
\end{align}
where the TT-cores are interpreted as third order tensors $(G_k)_{j,i,\ell}=G_k(i)_{j,\ell}$.

\begin{figure}[!h]
\centering
\includegraphics{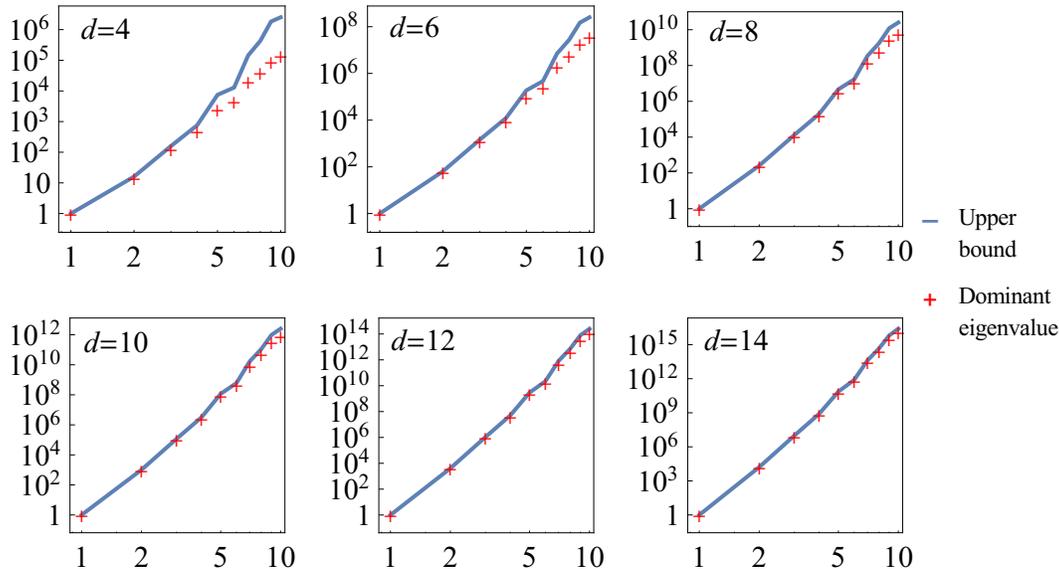}
\caption{The dominant eigenvalues of the LCM tensor $[S_d]$, $S=\{1,2,\ldots,n\}$, accompanied by the theoretical upper bound given by Proposition~\ref{evpbound} with increasing dimension $1\leq n\leq 10$ and for fixed order $d\in\{4,6,8,10,12,14\}$.}\label{eigfig}
\end{figure}

The dominant eigenvalues of the LCM tensor $[S_d]$, $S=\{1,\ldots,n\}$, were computed for $d\in\{4,6,8,10,12,14\}$ with $1\leq n\leq 10$ using Algorithm~\ref{power} and they are displayed in Figure~\ref{eigfig}. The tensor $[S_d]$ was represented using the TT-decomposition given by Theorem~\ref{ttdecomp} and the formulae \eqref{contr1} and \eqref{contr2} were used to compute the mode contractions in Algorithm~\ref{power} efficiently. The sets $S^{\vee k}$ were constructed for $1\leq k\leq d$ to determine the upper bound of Proposition~\ref{evpbound}. Algorithm~\ref{power} was run for $10\,000$ iterations, which led to small enough errors $\overline{\lambda_k}-\underline{\lambda_k}<\varepsilon$ to fall below the theoretical upper bound. The error bounds are not visible in the scale of Figure~\ref{eigfig} and are thus omitted.

The dominant eigenvalues of the LCM tensor increase at an exponential rate. This is a consequence of the fact that the values of the elements of $A=[S_d]$ increase exponentially with respect to $d$, which makes the operation $Ax^{d-1}$ act as a magnifier on the vector $x$. At the same time, the upper bound of Proposition~\ref{evpbound} scales with the action of tensor $A$ on vector $x$ due to the exponential term $n^{d-1}-1$ in Proposition~\ref{evpbound}, which dominates the growth rate of coefficients $c_k$.

\section{Conclusions}\label{concl}
In recent literature, the structure and properties of meet and join tensors have gained attention within the lattice-theoretic community. However, the properties of these tensors have not been discussed when the domain of definition is not closed with respect to the meet or join operation, respectively. In this paper, decompositions of general join tensors in both polyadic and tensor-train formats have been derived. In addition, it can be shown that both the CP-rank and TT-rank of general join tensors are bounded from below and the obtained decompositions are rank decompositions when the dimension $n$ is bounded by the order $d$ by $n\leq \lfloor d/2\rfloor$. With increasing order, the representation of join tensors---numerical or explicit---with dimensionality $n\leq \lfloor d/2\rfloor$ thus cannot be divorced from the cardinalities of these sets without imposing additional assumptions on the semilattice.

Decompositions of meet and join matrices have been an integral tool in the study of their properties, and this trend is likely to continue in the tensor case as well. In this paper, dominant eigenvalues of a class of LCM tensors were computed as an application of the obtained decompositions. Decompositions of join tensors may be used in the future to improve bounds on their eigenvalues and other characteristics, and we hope that the results of this paper instigate this discussion.

Low-rank tensor decompositions are a subject of intense study of the tensor community due to the intricacy of determining explicit rank decompositions of tensors. Since this work provides insight into the construction of tensor rank decompositions for a certain class of tensors, the results of this paper may be useful to this audience as well.

\subsection*{Acknowledgement}

The author is grateful for the invaluable advice given by Professor Pauliina Ilmonen during the writing of this article.

\subsection*{Disclosure statement}

No potential conflict of interest was reported by the author.

\subsection*{Funding}

The author has been supported by the Academy of Finland under Grant 267789.

\bibliographystyle{tfnlm}
\bibliography{jointensors}
\end{document}